\documentclass[a4paper,11pt]{article}
\usepackage[latin1]{inputenc}
\usepackage{amssymb,amsmath}
\usepackage{enumerate,enumitem}
\usepackage{graphicx}
\usepackage{color}
\usepackage{cite}
\usepackage{authblk}
\usepackage{tikz}

\newtheorem{theorem}{Theorem}[section]

\newtheorem{proposition}[theorem]{Proposition}

\newtheorem{corollary}[theorem]{Corollary}

\newtheorem{lemma}[theorem]{Lemma}

\newtheorem{question}[theorem]{Question}

\newcommand{\proof}{\noindent{\bf Proof.\ }}
\newcommand{\qed}{\hfill $\square$ \bigskip}

\newcommand{\comment}[1]{\textcolor{red}{#1}}

\begin{document}

\title{On the Djokovi\'c-Winkler relation and its closure in subdivisions of fullerenes, triangulations, and chordal graphs}

\author{
Sandi Klav\v zar $^{a,b,c}$\quad Kolja Knauer $^d$\quad Tilen Marc $^{a,c,e}$
}
\date{}
\maketitle

\vspace*{-1cm}
\begin{center}
$^a$ Faculty of Mathematics and Physics, University of Ljubljana, Slovenia \\
\medskip

$^b$ Faculty of Natural Sciences and Mathematics, University of Maribor, Slovenia \\
\medskip

$^{c}$ Institute of Mathematics, Physics and Mechanics, Ljubljana, Slovenia \\
\medskip

$^{d}$ Aix Marseille Univ, Universit\'e de Toulon, CNRS, LIS, Marseille, France\\
\medskip

$^{e}$ XLAB d.o.o., Ljubljana, Slovenia
\end{center}

\begin{abstract}
It was recently pointed out that certain SiO$_2$ layer structures and SiO$_2$ nanotubes can be described as full subdivisions aka subdivision graphs of partial cubes. A key tool for analyzing distance-based topological indices in molecular graphs is the Djokovi\'c-Winkler relation $\Theta$ and its transitive closure $\Theta^\ast$. In this paper we study the behavior of $\Theta$ and $\Theta^\ast$ with respect to full subdivisions. We apply our results to describe $\Theta^\ast$ in full subdivisions of fullerenes, plane triangulations, and chordal graphs.  
\end{abstract}

\noindent
{\bf E-mails}: sandi.klavzar@fmf.uni-lj.si, kolja.knauer@lis-lab.fr, tilen.marc@fmf.uni-lj.si


\medskip\noindent
{\bf AMS Subj. Class.}: 05C12, 92E10

\section{Introduction}
\label{sec:intro}

Partial cubes, that is, graphs that admit isometric embeddings into hypercubes, are of great interest in metric graph theory. Fundamental results on partial cubes are due to Chepoi~\cite{chepoi-1988}, Djokovi\'c~\cite{djokovic-1973}, and Winkler~\cite{winkler-1984}. The original source for their interest however goes back to the paper of Graham and Pollak~\cite{graham-1971}. For additional information on partial cubes we refer to the books~\cite{deza-1997, eppstein-2008}, the semi-survey~\cite{ovchinnikov-2008}, recent papers~\cite{albenque-2016, cardinal-2015, marc-2017}, as well as references therein.  

Partial cubes offer many applications, ranging from the original one in interconnection networks~\cite{graham-1971} to media theory~\cite{eppstein-2008}. Our motivation though comes from mathematical chemistry where many important classes of chemical graphs are partial cubes. In the seminal paper~\cite{klavzar-1995} it was shown that the celebrated Wiener index of a partial cube can be obtained without actually computing the distance between all pairs of vertices. A decade later it was proved in~\cite{klavzar-2006}, based on the Graham-Winkler's canonical metric embedding~\cite{graham-1985}, that the method extends to arbitrary graphs. The paper~\cite{klavzar-1995} initiated the theory under the common name ``cut method,'' while~\cite{klavzar-2015} surveys the results on the method until 2015 with 97 papers in the bibliography. The cut method has been further developed afterwards, see~\cite{crepnjak-2017, tratnik-2017, tratnik-2018} for some recent results on it related to partial cubes. 

Now, in a series of papers~\cite{arockiaraj-2019+, arockiaraj-2019, arockiaraj-2019++} it was observed that certain SiO$_2$ layer structures and SiO$_2$ nanotubes that are of importance in chemistry can be described as the full subdivisions aka subdivision graphs of relatively simple partial cubes. (The paper~\cite{tian-2018} can serve as a possible starting point for the role of SiO$_2$ nanostructures in chemistry.) The key step of the cut-method for distance based (as well as some other) invariants is to understand and compute the relation $\Theta^\ast$. Therefore in~\cite{arockiaraj-2019} it was proved that the $\Theta^\ast$-classes of the full subdivision of a partial cube $G$ can be obtained from the $\Theta^\ast$-classes of $G$. Note that in a partial cube the latter coincide with the $\Theta$-classes. 

The above developments yield the following natural, general problem that intrigued us: Given a graph $G$ and its $\Theta^\ast$-classes, determine the $\Theta^\ast$-classes of the full subdivision of $G$. In this paper we study this problem and prove several general results that can be applied in cases such as in~\cite{arockiaraj-2019+, arockiaraj-2019, arockiaraj-2019++} in mathematical chemistry as well as elsewhere. In the next section we list known facts about the relations $\Theta$ and $\Theta^\ast$ as well as the distance function in full subdivisions needed in the rest of the paper. In Section~\ref{sec:general-properties}, general properties of the relations $\Theta$ and $\Theta^\ast$ in full subdivisions are derived. These properties are then applied in the subsequent sections. In the first of them, $\Theta^\ast$ is described for fullerenes (a central class of chemical graph theory, see e.g.~\cite{AKS16,SWA15}) and plane triangulations. In Section~\ref{sec:chordal} the same problem is solved for chordal graphs.

\section{Preliminaries}
\label{sec:prelim}

If $R$ is a relation, then $R^\ast$ denotes its transitive closure. The distance $d_G(x,y)$ between vertices $x$ and $y$ of a connected graph $G$ is the usual shortest path distance. If $x\in V(G)$ and $e = \{y,z\}\in E(G)$, then let
$$d_G(x, e) = \min\{d_G(x,y), d_G(x,z)\}\,.$$ 
Similarly, if $e = \{x,y\}\in E(G)$ and $f = \{u,v\}\in E(G)$, then we set 
$$d_G(e,f) = \min\{d_G(x,u), d_G(x,v), d_G(y,u), d_G(y,v)\}\,.$$ 
Note that the latter function does not yield a metric space because if $e$ and $f$ are adjacent edges then $d_G(e,f) = 0$. To get a metric space, one can define the distance between edges as the distance between the corresponding vertices in the line graph of $G$. But for our purposes the function $d_G(e,f)$ as defined is more suitable.   

Edges $e=\{x,y\}$ and $f=\{u,v\}$ of a graph $G$ are in relation $\Theta$, shortly $e\Theta f$, if $d_G(x,u) + d_G(y,v) \not= d_G(x,v) + d_G(y,u)$. If $G$ is bipartite, then the definition simplifies as follows. 

\begin{lemma}
\label{lem:Theta-in-bipartite}
If $e=\{x,y\}$ and $f=\{u,v\}$ are edges of a bipartite graph $G$ with $e\Theta f$, then the notation can be chosen such that $d_G(u, x) = d_G(v, y) = d_G(u, y) - 1 = d_G(v, x) - 1$. 
\end{lemma}

The relation $\Theta$ is reflexive and symmetric. Hence $\Theta^\ast$ is thus an equivalence, its classes are called $\Theta^\ast$-classes. Partial cubes are precisely those connected bipartite graph for which $\Theta = \Theta^\ast$ holds~\cite{winkler-1984}. In partial cubes we may thus speak of $\Theta$-classes instead of $\Theta^\ast$-classes. In the following lemma we collect properties of $\Theta$ to be implicitly or explicitly used later on. 

\begin{lemma}
\label{lem:properties-of-Theta}
\begin{enumerate}
\item[(i)] If $P$ is a shortest path in $G$, then no two distinct edges of $P$ are in relation $\Theta$. 
\item[(ii)] If $e$ and $f$ are edges from different blocks of a graph $G$, then $e$ is not in relation $\Theta$ with $f$. 
\item[(iii)] If $e$ and $f$ are edges of an isometric cycle $C$ of a bipartite graph $G$, then $e\Theta f$ if and only if $e$ and $f$ are antipodal edges of $C$. 
\item[(iv)] If $H$ is an isometric subgraph of a graph $G$, then $\Theta_{H}$ is the restriction of $\Theta_{G}$ to $H$. 
\end{enumerate} 
\end{lemma}

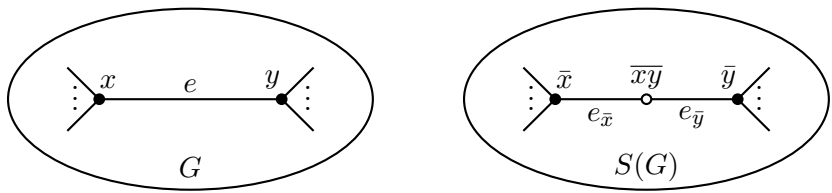
\begin{figure}[ht!]
\begin{center}
\begin{tikzpicture}[scale=0.6,style=thick,x=1cm,y=1cm]
\def\vr{3pt}

\coordinate (x) at (0,0);
\coordinate (y) at (4,0);
\coordinate (xbar) at (10,0);
\coordinate (ybar) at (14,0);
\coordinate (xybar) at (12,0);
\draw (x) -- (y);
\draw (xbar) -- (ybar);
\draw (x) -- ++(-0.7,-0.7)(x);
\draw (x) -- ++(-0.7,0.7)(x);
\draw (y) -- ++(0.7,-0.7)(y);
\draw (y) -- ++(0.7,0.7)(y);
\draw (xbar) -- ++(-0.7,-0.7)(xbar);
\draw (xbar) -- ++(-0.7,0.7)(xbar);
\draw (ybar) -- ++(0.7,-0.7)(ybar);
\draw (ybar) -- ++(0.7,0.7)(ybar);
\draw (2,0) ellipse (4cm and 2cm);
\draw (12,0) ellipse (4cm and 2cm);
\draw (x)[fill=black] circle(\vr);
\draw (y)[fill=black] circle(\vr);
\draw (xbar)[fill=black] circle(\vr);
\draw (ybar)[fill=black] circle(\vr);
\draw (xybar)[fill=white] circle(\vr);
\draw (x) ++(0.2,0) node[above] {$x$};
\draw (x) ++(-0.2,0.2) node[left] {$\vdots$};
\draw (y) ++(-0.2,0) node[above] {$y$};
\draw (y) ++(0.9,0.2) node[left] {$\vdots$};
\draw (2,0.3) node {$e$};
\draw (xbar) ++(0.2,0) node[above] {$\bar{x}$};
\draw (xbar) ++(-0.2,0.2) node[left] {$\vdots$};
\draw (ybar) ++(-0.2,0) node[above] {$\bar{y}$};
\draw (ybar) ++(0.8,0.2) node[left] {$\vdots$};
\draw (xybar) node[above] {$\overline{xy}$};
\draw (xybar) ++(-1,-0.4) node {$e_{\bar{x}}$};
\draw (xybar) ++(1,-0.4) node {$e_{\bar{y}}$};
\draw (2,-1.5) node {$G$};
\draw (12,-1.5) node {$S(G)$};
\end{tikzpicture}
\end{center}
\caption{Notation for the vertices and edges of $S(G)$.}
\label{fig:notation}
\end{figure}

If $G$ is a graph, then the graph obtained from $G$ by subdividing each each of $G$ exactly once is called the {\em full subdivision (graph)} of $G$ and denoted with $S(G)$.  We will use the following related notation. If $x\in V(G)$ and $e = \{x,y\}\in E(G)$, then the vertex of $S(G)$ corresponding to $x$ will be denoted by $\bar{x}$ and the vertex of $S(G)$ obtained by subdividing the edge $e$ with $\overline{xy}$. Two edges incident with $\overline{xy}$ will be denoted with $e_{\bar{x}}$ and  $e_{\bar{y}}$, where $e_{\bar{x}} = \{\bar{x}, \overline{xy}\}$ and $e_{\bar{y}} = \{\bar{y}, \overline{xy}\}$. See Fig.~\ref{fig:notation} for an illustration.

The following lemma is straightforward, cf.~\cite[Lemma 2.3]{klavzar-2003}. 
\begin{lemma}
\label{lem:distances}
If $G$ is a connected graph, then the following assertions hold.  
\begin{enumerate}
\item[(i)] If $x,y\in V(G)$, then $d_{S(G)}(\bar{x}, \bar{y}) = 2d_G(x,y)$. 
\item[(ii)] If $x\in V(G)$ and $\{y,z\}\in E(G)$, then $d_{S(G)}(\bar{x}, \overline{yz}) = 2 d_G(x,\{y,z\}) + 1$. 
\item[(iii)] If $\{x,y\}, \{u,v\}\in E(G)$, then $d_{S(G)}(\overline{xy}, \overline{uv}) = 2 d_G(\{x,y\}, \{u,v\}) + 2$. 
\end{enumerate}
\end{lemma}

\section{$\Theta^*$ in full subdivisions}
\label{sec:general-properties}

\begin{lemma}\label{lem:from-S(G)-to-G}
If $G$ is a connected graph and $e_{\bar{x}}\, \Theta_{S(G)}\, f_{\bar{u}}$, then $e\, \Theta_G\, f$. 
\end{lemma}

\proof
Let $e = \{x,y\}$ and $f = \{u, v\}$. If $\bar{x} = \bar{u}$ and $\bar{y} = \bar{v}$, then $e_{\bar{x}} = f_{\bar{u}}$ and $e = f$, so there is nothing to prove. If $\bar{x} = \bar{v}$ and $\bar{y} = \bar{u}$, then $e_{\bar{x}}$ and $f_{\bar{u}}$ are adjacent edges which cannot be in relation $\Theta_{S(G)}$ because $S(G)$ is triangle-free. For the same reason the situation $\bar{x} =  \bar{u}$ and $\bar{y} \ne \bar{v}$ is not possible. Assume next that $\bar{x} =  \bar{v}$ and $\bar{y} \ne \bar{u}$. Then $d_{S(G)}(\bar{u}, \overline{xy}) = 3$ by Lemma~\ref{lem:distances}, and hence $\overline{xy}, \bar{x}, \overline{uv}, \bar{u}$ is a geodesic containing $e_{\bar{x}}$ and $f_{\bar{u}}$, contradiction the assumption $e_{\bar{x}}\, \Theta_{S(G)}\, f_{\bar{u}}$. In the rest of the proof we may thus assume that  $\{x, y\}\cap \{u,v\} = \emptyset$. 

Since $S(G)$ is bipartite, in view of Lemma~\ref{lem:Theta-in-bipartite} we need to consider the following two cases, where, using Lemma~\ref{lem:distances}(i), we can assume that the distances $d_{S(G)}(\bar{x}, \bar{u})$ and $d_{S(G)}(\overline{xy}, \overline{uv})$ are even. Based on the assumption $e_{\bar{x}}\, \Theta_{S(G)}\, f_{\bar{u}}$, we have $d_{S(G)}(\bar{x}, \bar{u}) + d_{S(G)}(\overline{xy}, \overline{uv}) = d_{S(G)}(\bar{x}, \overline{uv}) + d_{S(G)}(\overline{xy}, \bar{u})$ in a bipartite graph, thus the following cases.

\medskip\noindent
{\bf Case 1.} $d_{S(G)}(\bar{x}, \bar{u}) = d_{S(G)}(\overline{xy}, \overline{uv}) = 2k$ and $d_{S(G)}(\bar{x}, \overline{uv}) = d_{S(G)}(\overline{xy}, \bar{u}) = 2k+1$.   \\
In the following, Lemma~\ref{lem:distances} will be used all the time.

By $2k = d_{S(G)}(\overline{xy}, \overline{uv}) = 2 d_G(\{x,y\}, \{u,v\}) + 2$, we get
$$k-1\leq d_G(y,v), d_G(x,u),d_G(x,v),d_G(y,u),$$
where the lower bound is attained at least once.

Since $d_{S(G)}(\bar{x}, \bar{u}) = 2k$, we have $d_G(x,u) = k$.
Because $d_{S(G)}(\bar{x}, \overline{uv}) = 2k+1$, we find that $d_G(x,\{u,v\}) = k$ and hence in particular $d_G(x,v) \ge k$. Similarly, as $ d_{S(G)}(\overline{xy}, \bar{u}) = 2k+1$ we have $d_G(u,\{x,y\}) = k$ and hence in particular $d_G(u,y) \ge k$. With the first observation this yields $k-1=d_G(y,v)$.
In summary, 
$$d_G(x,u) + d_G(y,v)= k + (k-1) \not= k + k \le d_G(x,v) + d_G(y,u)\,,$$
which means that $e\, \Theta_G\, f$. 

\medskip\noindent
{\bf Case 2.} $d_{S(G)}(\bar{x}, \bar{u}) = d_{S(G)}(\overline{xy}, \overline{uv}) = 2k$ and $d_{S(G)}(\bar{x}, \overline{uv}) = d_{S(G)}(\overline{xy}, \bar{u}) = 2k-1$.   \\
Again, $d_{S(G)}(\bar{x}, \bar{u}) = 2k$ implies $d_G(x,u) = k$. The assumption $d_{S(G)}(\bar{x}, \overline{uv}) = 2k-1$ yields $d_G(x,\{u,v\}) = k-1$ and consequently $d_G(x,v) = k-1$. The condition $d_{S(G)}(\overline{xy}, \bar{u}) = 2k-1$ implies $d_G(u,\{x,y\}) = k-1$ and so $d_G(u,y) = k-1$. Finally, the assumption $d_{S(G)}(\overline{xy}, \overline{uv}) = 2k$ gives us $d_G(\{x,y\}, \{u,v\}) = k-1$, in particular, $d_G(y,v) \ge  k-1$. Putting these facts together we get 
$$d_G(x,u) + d_G(y,v) \ge k + (k-1) > (k-1) + (k-1) = d_G(x,v) + d_G(y,u)\,,$$
hence again $e\, \Theta_G\, f$. 
\qed

Lemma~\ref{lem:from-S(G)-to-G} implies the following result on the relation $\Theta^\ast$. 

\begin{corollary}
\label{cor:refinement}
If $e_{\bar{x}}\,\Theta^\ast_{S(G)}\,f_{\bar{u}}$, then $e\Theta_G^\ast f$. 
\end{corollary}

\proof 
Suppose $e_{\bar{x}}\,\Theta^\ast_{S(G)}\,f_{\bar{u}}$. Then there exists a positive integer $k$ such that 
$$e_{\bar{x}}\,\Theta_{S(G)}\,f^{(1)}_{\overline{x}_1}, 
f^{(1)}_{\overline{x_1}}\,\Theta_{S(G)}\,f^{(2)}_{\overline{x}_2},
\ldots,  
f^{(k)}_{\overline{x}_k}\,\Theta_{S(G)}\,f_{\bar{u}}\,.$$ 
Then, by Lemma~\ref{lem:from-S(G)-to-G}, we have 
$$e\,\Theta_{G}\,f^{(1)}, 
f^{(1)}\,\Theta_{G}\,f^{(2)},
\ldots,  
f^{(k)}\,\Theta_{G}\,f\,,$$ 
implying that $e\Theta_G^\ast f$.
\qed

The next lemma is a partial converse to Lemma~\ref{lem:from-S(G)-to-G}. 

\begin{lemma}\label{lem:converserefinement}
If $e\,\Theta_G\, f$, then there is a pair of edges $e_{\bar{x}}, f_{\bar{u}}$ in $S(G)$ such that $e_{\bar{x}}\,\Theta_{S(G)}\,f_{\bar{u}}$. Moreover, if $G$ is bipartite, then there are two (disjoint) such pairs.
\end{lemma}

\proof
Let $e = \{x,y\}$, $f = \{u, v\}$, and let $k = d_G(x,u)$. Since $e\,\Theta_G\, f$, we may without loss of generality assume that $d_G(x,u) + d_G(y,v) < d_G(y,u) + d_G(x,v)$ and that $d_G(x,u) \le d_G(y,v)$. We distinguish the following cases. 

\medskip\noindent
{\bf Case 1.} $d_G(y,v) = k$.\\
In this case, $\{d_G(x,v), d_G(y,u)\} \subseteq \{k-1,k,k+1\}$. Moreover, our assumption about the sum of distances implies that $\{d_G(x,v), d_G(y,u)\} \subseteq \{k,k+1\}$. Since $e\,\Theta_G\, f$, the two distances cannot both be equal to $k$. Hence, up to symmetry, we need to consider the following two subcases.  

Suppose $d_G(x,v) = d_G(y,u) = k+1$. Then $d_{S(G)}(\bar{x}, \bar{v}) = 2k + 2$, $d_{S(G)}(\overline{xy}, \overline{uv}) = 2k +2$, $d_{S(G)}(\bar{x}, \overline{uv}) = 2k + 1$, and $d_{S(G)}(\overline{xy}, \bar{v}) = 2k+1$. Hence  $e_{\bar{x}}\,\Theta_{S(G)}\,f_{\bar{v}}$.

Suppose $d_G(x,v) = k$ and $d_G(y,u) = k + 1$. Then $d_{S(G)}(\bar{y}, \bar{u}) = 2k + 2$, $d_{S(G)}(\overline{xy}, \overline{uv}) = 2k +2$, $d_{S(G)}(\bar{y}, \overline{uv}) = 2k + 1$, and $d_{S(G)}(\overline{xy}, \bar{u}) = 2k+1$. Hence  $e_{\bar{y}}\,\Theta_{S(G)}\,f_{\bar{u}}$. A similar situation occurs when $d_G(x,v) = k + 1$ and $d_G(y,u) = k$.

\medskip\noindent
{\bf Case 2.} $d_G(y,v) = k + 1$.\\
Again, $\{d_G(x,v), d_G(y,u)\} \subseteq \{k-1,k,k+1\}$, but since $d_G(x,u) + d_G(y,v) < d_G(y,u) + d_G(x,v)$ it must be that $d_G(x,v) = d_G(y,u) = k+1$. Then $d_{S(G)}(\bar{y}, \bar{v}) = 2k + 2$, $d_{S(G)}(\overline{xy}, \overline{uv}) = 2k +2$, $d_{S(G)}(\bar{y}, \overline{uv}) = 2k + 3$, and $d_{S(G)}(\overline{xy}, \bar{v}) = 2k+3$. Hence  $e_{\bar{y}}\,\Theta_{S(G)}\,f_{\bar{v}}$.

\medskip\noindent
{\bf Case 3.} $d_G(y,v) = k + 2$.\\
In this case the fact that $\{d_G(x,v), d_G(y,u)\} \subseteq \{k-1,k,k+1\}$ implies that $d_G(x,u) + d_G(y,v) \geq d_G(y,u) + d_G(x,v)$. As this is not possible, the first assertion of the lemma is proved. 

\medskip
Assume now that $G$ is bipartite. Combining Lemma~\ref{lem:Theta-in-bipartite} with the above case analysis we infer that the only case to consider is when $d_G(x,u) = d_G(y,v) = k$ and $d_G(x,v) = d_G(y,u) = k+1$. Then, just in the first subcase of the above Case 1 we get that 
$e_{\bar{x}}\,\Theta^\ast_{S(G)}\,f_{\bar{v}}$ and, similarly, 
$e_{\bar{y}}\,\Theta^\ast_{S(G)}\,f_{\bar{u}}$. 
\qed

We say that cycles $C$ and $C'$ of $G$ are \emph{isometrically touching} if $|E(C)\cap E(C')| = 1$  and $C\cup C'$ is an isometric subgraph of $G$. Note that isometrically touching cycles are isometric.

\begin{figure}[ht!]
\centering
\includegraphics[width=.6\textwidth]{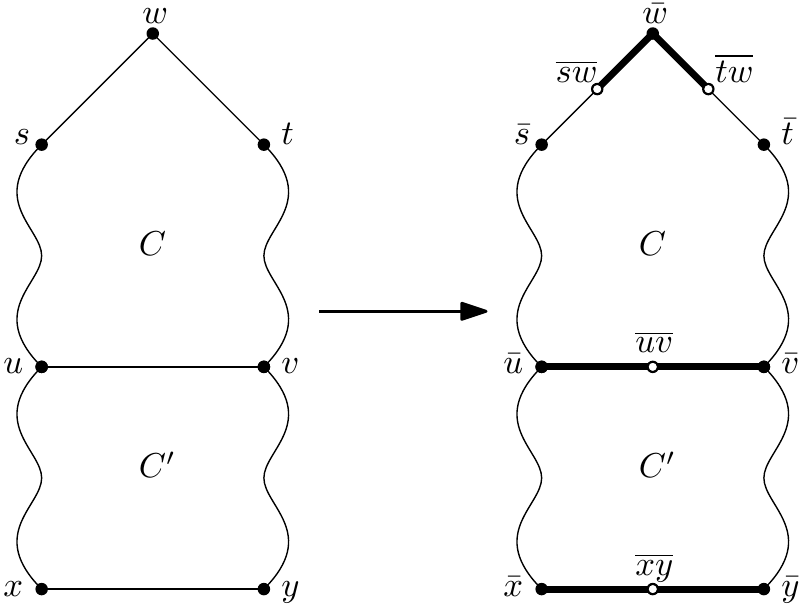}
\caption{Isometrically touching cycles and their subdivisions.}
\label{fig:cycles}
\end{figure}

\begin{lemma}\label{lem:cycles}
Let $C$ and $C'$ be isometrically touching cycles in $G$ with $E(C)\cap E(C') = \{e\}$. Then in $S(G)$ both edges corresponding to $e$ are in the same $\Theta^*_{S(G)}$-class. Moreover, this class contains the edges thickened in Fig.~\ref{fig:cycles}.
\end{lemma}

\proof
 We take the notation from Fig.~\ref{fig:cycles} and content ourselves with only providing the proof for the case where $C$ is odd and $C'$ is even. The other cases go through similarly. From Lemma~\ref{lem:properties-of-Theta}(iii) we get that $\{\bar{u},\overline{uv}\}\Theta_{S(G)}\{\bar{w},\overline{tw}\}$ and $\{\bar{u},\overline{uv}\}\Theta_{S(G)}\{\bar{y},\overline{xy}\}$. However, note now that $d(\bar{y},\bar{w})=d(\overline{xy},\overline{sw})=d(\bar{y},\overline{sw})-1=d(\overline{xy},\bar{w})-1$. Thus we also have $\{\bar{w},\overline{sw}\}\Theta_{S(G)}\{\bar{y},\overline{xy}\}$. Since $\{\bar{w},\overline{sw}\}$ is also in relation with $\{\bar{v},\overline{uv}\}$ we obtain the claim for $\Theta^*_{S(G)}$ by taking the transitive closure.
\qed 

For the full subdivision $S(G)$ of $G$ denote by $S(\Theta^*_G)$, the relation on the edges of $S(G)$, where $\{\bar{x},\overline{xy}\}$ and $\{\bar{u},\overline{uv}\}$ are in relation $S(\Theta^*_G)$ if and only if $\{x,y\}\Theta^*\{u,v\}$. In particular, $\{\bar{x},\overline{xy}\}$ and $\{\overline{xy},\bar{y}\}$ are always in relation. 

\begin{lemma}\label{lem:halfedgesthesame}
 We have $\{\bar{x},\overline{xy}\}\Theta^*_{S(G)}\{\overline{xy},\bar{y}\}$ for all $\{x,y\}\in G$ if and only if $\Theta^*_{S(G)}=S(\Theta^*_G)$.
\end{lemma}

\proof
 The backwards direction holds by definition. Conversely, by Lemma~\ref{lem:from-S(G)-to-G} we have that if $\{\bar{x},\overline{xy}\}\Theta^*_{S(G)}\{\overline{uv},\bar{v}\}$, 
 then $\{x,y\}\Theta^*\{u,v\}$. Therefore, $\Theta^*_{S(G)}\subseteq S(\Theta^*_G)$. On the other hand, Lemma~\ref{lem:converserefinement} assures that if $\{x,y\}\Theta^*\{u,v\}$, then there is a pair 
 $\{\bar{x},\overline{xy}\}
 \Theta^*_{S(G)}\{\overline{uv},\bar{v}\}$, 
 but then by our assumption also $\{\bar{y},\overline{xy}\}\Theta^*_{S(G)}\{\overline{uv},\bar{v}\}$ and so on. Thus, $\Theta^*_{S(G)}\supseteq S(\Theta^*_G)$.
\qed 

Lemma~\ref{lem:cycles} and~\ref{lem:halfedgesthesame} immediately yield:
\begin{proposition}\label{prop:together}
 If every edge of $G$ is in the intersection of two isometrically touching cycles, then $\Theta^*_{S(G)}=S(\Theta^*_G)$.
\end{proposition}

\section{$\Theta^\ast$ in subdivisions of fullerenes and plane triangulations}
\label{sec:fullerenes-triangulations}

In this section we study relation $\Theta^\ast$ in full subdivisions of fullerenes and plane triangulations, for which Proposition~\ref{prop:together} will be essential. We begin with fullerenes. Recall that a \emph{fullerene} is a cubic planar graph all of whose faces are of length $5$ or $6$.

A cycle $C$ of a connected graph $G$ is {\em separating} if $G\setminus C$ is disconnected and that a {\em cyclic edge-cut} of $G$ is an edge set $F$ such that $G\setminus F$ separates two cycles. To prove our main result on fullerenes we need the following result.

\begin{figure}[ht!]
\centering
\includegraphics[width=.6\textwidth]{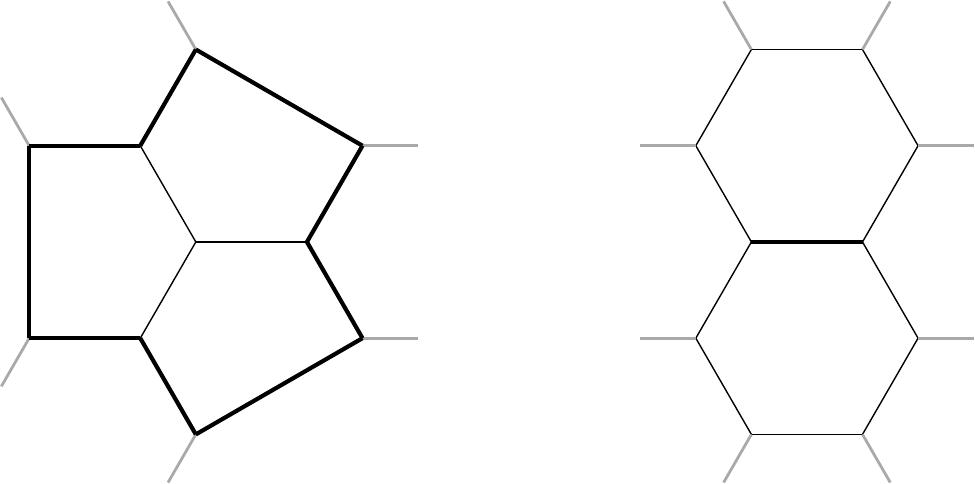}
\caption{A separating $9$-cycle and two isometrically touching $6$-cycles in a fullerene.}
\label{fig:C9}
\end{figure}

\begin{lemma}\label{lem:separating}
Given a fullerene graph $G$, every separating cycle of $G$ is  of length at least 9. Moreover, the only separating cycles of length 9, are the cycles separating a vertex incident only to 5-faces, see the left of Fig.~\ref{fig:C9}.
\end{lemma}

\proof
Let $C$ be separating cycle of length at most 9. Since $G$ is cubic, there are $|C|$ edges of $G$ incident to $C$ which are not in $C$. Thus, without loss of generality, we may assume at most four of them are in the inner side of $C$. As they form an edge cut, and since fullerenes are cyclically 5-edge-connected~\cite{doslic-2003}, the subgraph induced by vertices in the inner part of $C$ is a forest, say $F$. If $F$ consists of one vertex, say $v$, then we have three edges connecting $v$ to $C$ which form three faces $G$. As each of these faces is of length at least 5, they are exactly 5-faces. Otherwise, $F$ contains at least two vertices $u$ and $v$ each of which is either an isolated vertex of $F$ or a leaf. As they are of degree 3 in $F$, each of them must be connected by two edges to $C$. And since there at most four such edges, it follows that $u$ and $v$ are of degree 1 in $F$ and that every other vertex of $F$ is of degree 3 in $F$, which means there no other vertex and $u$ is adjacent to $v$. Thus inside $C$ we have five edges and four faces. But $C$ itself is of length at most 9 and thus one of these four faces is of length at most 4, a contradiction with the choice of $G$.
\qed

%
%
%
%
\begin{theorem}
If $G$ is a fullerene, then $\Theta^*_{S(G)}=S(\Theta^*_G)$.
\end{theorem}

\proof
We claim that every edge $e$ of $G$ is the intersection of two isometrically touching cycles. For this sake consider the cycles $C$ and $C'$ that lie on the boundary of the faces containing $e$. We have to prove that the union $C\cup C'$ is isometric. Assume on the contrary that this is not the case, that is, there exist vertices $u,v\in C\cup C'$ such that there is a shortest $u,v$-path $P$ (in $G$) interiorly disjoint from $C\cup C'$, that is shorter than any shortest path $P'$ from $u$ to $v$ in $C\cup C'$. 

Consider the cycle $C''$ obtained by joining $P$ and a shortest path $P'$ from $u$ to $v$ in $C\cup C'$. Since $C$ and $C'$ are of length at most $6$, the graph $C\cup C'$ is of diameter at most $5$, thus the cycle $C''$ is of length at most $9$. 
We will prove that there is a separating cycle contradicting Lemma~\ref{lem:separating}.

First, note that if $e\in P'$, then $C''$ separates the graph $C\cup C'$. Thus, by Lemma~\ref{lem:separating} $C''$ is of length at most $9$, so the endpoints of $P'$ must be at distance $5$ on $C\cup C'$, i.e., one is in $C$ and the other in $C'$. Thus, both sides of $C''$ contain more than one vertex, contradicting Lemma~\ref{lem:separating}.

Hence, $P'$ is on the boundary of $C\cup C'$. Suppose that $C''$ is not induced. Then since the girth of fullerenes is $5$, there is a single chord from $P$ to $P'$ which splits $C''$ into a $5$-cycle $A$ and into a $5$- or a $6$-cycle $B$. In particular, $|C''|\geq 8$ and $P'$ has at least five vertices on $C''$. Thus, one vertex of $P'$ has degree $2$ in $C\cup C'$ and is not incident to the chord. Thus, this vertex has a neighbor in the interior of $A$ or $B$, that is, one of them is separating, contradicting Lemma~\ref{lem:separating}.

If $C''$ is induced, it follows from the fact that $C$ and $C'$ are faces and $|C''|\geq 5$, that $C''$ is not a face, i.e., it is separating. Thus, $|C''|=9$ and the patch $Q$ consisting of $C''$ and its interior is a single vertex surrounded by three 5-faces, see Lemma~\ref{lem:separating}. Moreover, $C$ and $C'$ are 6-faces so that their union can have diameter $5$. Note that any path $P'$ on the boundary of $C\cup C'$ of length $5$ uses only one vertex of degree $3$, see the right of Fig.~\ref{fig:C9}. But any path $P'$ of length $5$ on the boundary $Q$ uses at least two vertices of degree $2$, see the left of Fig.~\ref{fig:C9}. Thus, $P'$ cannot be in both boundaries simultaneously -- contradiction.
 
We have shown the claim from the beginning and Proposition~\ref{prop:together} yields the result.
\qed 

\begin{figure}[ht!]
\centering
\includegraphics[width=.6\textwidth]{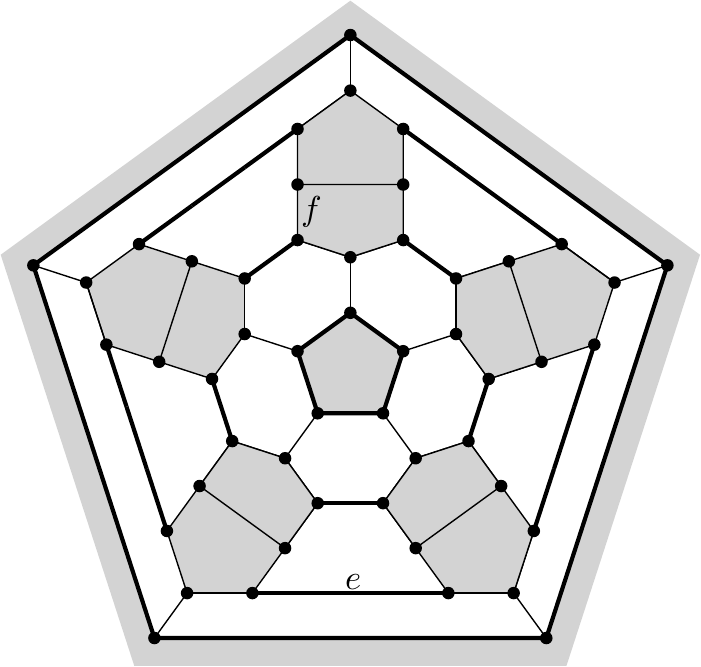}
\caption{A fullerene $G$ which has two $\overline{\Phi}^{\,\ast}$-classes (bold and normal edges), but only one $\Theta^*_G$-class since $e\Theta_G f$.}
\label{fig:counterxmpl}
\end{figure}

We have proved how $\Theta^*_G$ of a fullerene behaves with respect to subdivision. What can we say about $\Theta^*_G$ itself?
If $G$ is a fullerene, then we define a relation $\Phi$ on $E(G)$ as follows: $e\Phi f$ if $e$ and $f$ are opposite edges of a facial $C_6$. Relation $\Phi$ falls into cycles and paths, that have been called \emph{railroads}~\cite{deza-2004}. In particular, it has been shown that cycles can have multiple self-intersection. We denote by $\overline{\Phi}$ the relation where additionally any two non-incident edges of a facial $C_5$ are in relation. Finally, recall that $\overline{\Phi}^{\,\ast}$ denotes the transitive closure of $\overline{\Phi}$. Since faces are isometric subgraphs, it is easy to see that $\overline{\Phi}$ is a refinement of $\Theta_G$ as well as $\overline{\Phi}^{\,\ast}$ is a refinement of $\Theta^*_G$. One might believe that the converse also holds, but the example in Fig.~\ref{fig:counterxmpl} shows that this is not always the case. We believe that determining $\Theta^*_G$ in fullerenes is an interesting problem.

\bigskip

We now turn our attention to plane triangulations. It is straightforward to verify that if $G$ is a plane triangulation, then $\Theta^*$ consists of a single class. On the other hand, $\Theta^*$ on the full subdivision of a plane triangulation has the following non-trivial structure. 

\begin{theorem}
 Let $G\neq K_4$ be a plane triangulation. Then $\Theta^*_{S(G)}$ consists of one global class $\gamma$, plus one class $\gamma_x$ for every degree three vertex $x$. Here, if $N(x) = \{y_1, y_2, y_3\}$, then $\gamma_x = \{\{\bar{y}_1,\overline{y_1x}\},\{\bar{y}_2,\overline{y_2x}\},\{\bar{y}_3,\overline{y_3x}\}\}$. If $G=K_4$ the same holds, except that there is no global class $\gamma$.
\end{theorem}

\begin{figure}[ht!]
\centering
\includegraphics[width=.6\textwidth]{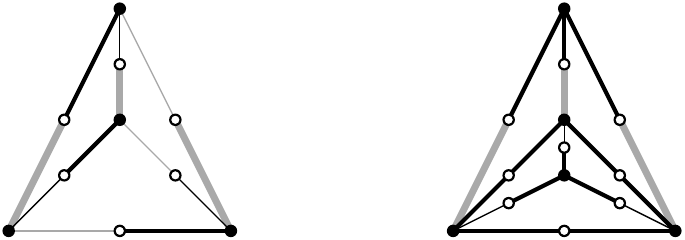}
\caption{The relation $\Theta^*$ in $S(K_4)$ and the full division of the graph obtained by stacking into one face.}
\label{fig:K4}
\end{figure}

\proof
Recall that $S(K_4)$ is a partial cube, cf.~\cite{klavzar-2003}, its $\Theta$-classes (= $\Theta^\ast$-classes) are shown In Fig.~\ref{fig:K4}. Hence the result holds for $K_4$. 

 We proceed by induction on the number of vertices. Let $G$ have minimum degree at least $4$, and let $e=\{x,y\}$ be an edge shared by triangles $C$ and $C'$ bonding faces of $G$. If  $C\cup C'$ is isometric, then by Lemma~\ref{lem:cycles} we have $\{\bar{x},\overline{xy}\}\Theta^*_{S(G)}\{\overline{xy},\bar{y}\}$. Otherwise, $C\cup C'$ induces a $K_4$, but since the minimum degree of $G$ is at least $4$, the other two triangles of the $K_4$ cannot be faces. An easy application of Lemma~\ref{lem:cycles} on the other edges of this $K_4$ implies $\{\bar{x},\overline{xy}\}\Theta^*_{S(G)}\{\overline{xy},\bar{y}\}$. Since in a triangulation there is only one $\Theta^*$-class, Proposition~\ref{prop:together} implies the result, that is, there is only one global class $\gamma$ in $S(G)$. 
 
 Now suppose that $G$ contains a vertex $v$ of degree $3$. The graph $G'=G\setminus \{v\}$ is a plane triangulation, thus our claim holds for $G'$ by induction. In particular, if $G'=K_4$, see Fig.~\ref{fig:K4} again. Otherwise, since $S(G')$ is an isometric subgraph of $S(G)$, Lemma~\ref{lem:properties-of-Theta}(iv) says that $\Theta_{S(G')}$ is the restriction of $\Theta_{S(G)}$ to $S(G')$. 
 
 Consider an edge $e=\{x,y\}$ of the triangle of $G$ that contains $v$. Note that the facial triangles $C,C'$ containing $e$ have an isometric union, so by Lemma~\ref{lem:cycles} we have $\{\bar{x},\overline{xy}\}\Theta^*_{S(G)}\{\overline{xy},\bar{y}\}$, which corresponds to our claim, since neither $x$ or $y$ can be of degree $3$. If one of them---say $x$---was of degree $3$ in $G'$, then now only the class $\gamma_x$ and $\gamma$ where merged. Since $G'\neq K_4$, not both $x$ and $y$ are of degree $3$. Note furthermore that by Lemma~\ref{lem:cycles} the edges incident to $v$ will all be in the class $\gamma$. 
 
 Finally, all the edges of the form $f=\{\bar{x},\overline{vx}\}$ are in relation $\Theta$ with each other. In order to see that they are the only constituents of the class $\gamma_v$ it suffices to notice that $d(\overline{vx}, z)=d(\overline{x}, z)+1$ for all $z\in S(G')$. The result then follows by Lemma~\ref{lem:properties-of-Theta} (i).
\qed

\section{$\Theta^\ast$ in subdivisions of chordal graphs}
\label{sec:chordal}

%

Recall that a graph is \emph{chordal} if all its induced cycles are of length $3$.
Similarly as in fullerenes we shall define relation $\Phi$ on the edges of $S(G)$, by $e\Phi f$ if $e,f$ are opposite edges of a $C_6$.

\begin{lemma}\label{lem:local}
If  $G$ is a chordal graph, then $\Phi_{S(G)}^* = \Theta_{S(G)}^*$.
\end{lemma}

\proof
Let $e \Theta_{S(G)} f$, where $e$ and $f$ are edges created by subdividing $\{a, b\}, \{c, d\} \in E(G)$, respectively. Then by Lemma \ref{lem:from-S(G)-to-G} we have $\{a, b\} \Theta \{c, d\}$. Similarly as in the proof of Lemma \ref{lem:converserefinement}, we have (up to symmetry) two options. 

\medskip\noindent
{\bf Case 1.} $d_{G}(a,c) = d_{G}(b, d)=k$.\\
We can assume that $d_{G}(a,d) \in \{k, k+1\}$ and $d_{G}(b, c) = k+1$. Let $P=p_0p_1\ldots p_k$ and $P'=p_0'p_1'\ldots p_k'$ be shortest $a,c$- and $b,d$-paths, respectively. Clearly, $P$ and $P'$ must be disjoint since otherwise it cannot hold $d_{G}(a,d) \in \{k, k+1\},  d_{G}(b, c) = k+1$. The cycle $C$ formed by $\{a,b\}, P', \{d,c\}, P$ must have a chord. Inductively adding chords we can show that there is a chord of $C$ incident with $a$ or $b$. Since $P$ and $P'$ are shortest paths and the assumptions on distances hold, it follows that the latter chord must be incident with $a$ and the vertex $p_1'$ of $P'$. In particular, $d_{G}(a,d) = k$. Similarly, one can show that there must be a chord between $p_1'$ and $p_1$, and inductively between every $p_ip_{i+1}'$ for $0\leq i < k$ and every $p_{i+1}p_{i+1}'$ for $0\leq i < k - 1$.

By the assumption on the distances, the only pair of subdivided edges of $\{a, b\}, \{c, d\}$, that is in relation $\Theta_{S(G)}$, is $\{\overline{b},\overline{ba}\} \Theta_{S(G)} \{\overline{c},\overline{cd}\}$, i.e., $e= \{\overline{b},\overline{ba}\}$ and $f= \{\overline{c},\overline{cd}\}$. Then 
$$\{\overline{b},\overline{ba}\} \Phi_{S(G)} \{\overline{a},\overline{ap_1'}\}\Phi_{S(G)} \{\overline{p_1},\overline{p_1p_1'}\} \Phi_{S(G)}\ \cdots\ \Phi_{S(G)} \{\overline{c},\overline{cd}\}\,.$$

\medskip\noindent
{\bf Case 2.} $d_{G}(a,c) = k, d_{G}(b, d)=k+1$.\\
Then we have $d_{G}(a,d) =  d_{G}(b, c) = k+1$. Similarly as above, shortest $a,c$- and $b,d$-paths, say $P=p_0p_1\ldots p_k$ and $P'=p_0'p_1'\ldots p_{k+1}'$, cannot intersect. Using the same notation as above, $C$ must have a chord incident with $a$ or $b$. By similar arguments, there must be a chord between every $p_ip_{i+1}'$ and $p_{i+1}p_{i+1}'$ for $0\leq i < k$.

By the assumption on the distances, the only pair of subdivided edges of $\{a, b\}, \{c, d\}$, that is in relation $\Theta_{S(G)}$, is $\{\overline{b},\overline{ba}\} \Theta_{S(G)} \{\overline{d},\overline{dc}\}$, i.e., $e= \{\overline{b},\overline{ba}\}$ and $f= \{\overline{d},\overline{dc}\}$. Then 
$$\{\overline{b},\overline{ba}\} \Phi_{S(G)} \{\overline{a},\overline{ap_1'}\}\Phi_{S(G)} \{\overline{p_1},\overline{p_1p_1'}\} \Phi_{S(G)}\ \cdots\ \Phi_{S(G)} \{\overline{d},\overline{dc}\}\,.$$
We have proved that $\Theta_{S(G)}\subset \Phi_{S(G)}^*$, thus $\Theta_{S(G)}^*= \Phi_{S(G)}^*$.
\qed

An edge of a graph $G$ is called {\em exposed} if it is properly contained in a single maximal complete subgraph of $G$. (This concept was recently introduced in~\cite{culbertson-2019+}, where it was proved that a $G$ is a connected chordal graph if and only if $G$ can be obtained from a complete graph by a sequence of removal of exposed edges.)  
%
Denote by $G^{-ee}$, for a chordal graph $G$, the graph obtained from $G$ by removing all its exposed edges. We will denote by $\textrm{c}(G^{-ee})$ the number of connected components of $G^{-ee}$. Note that the singletons of $G^{-ee}$ include the simplicial vertices of $G$, and if $G$ is 2-connected, its simplicial vertices coincide with singletons of $G^{-ee}$. 
It is straightforward to verify that if $G$ is a chordal graph, then $\Theta^*$ consists of a single class.  On the other hand, $\Theta^*$ on the full subdivision of a chordal graph has the following non-trivial structure. 

\begin{theorem}
Let $G$ be a $2$-connected, chordal graph. Then the coloring, that for an edge $\{a,b\}$ with $a$ being in the $i$-th connected component of $G^{-ee}$ colors edge $\{\overline{ab}, b\}$ with color $i$, corresponds to the $\Theta^*_{S(G)}$-partition. In particular, $|\Theta^*_{S(G)}| = \textrm{c}(G^{-ee})$.
\end{theorem}

\proof
We first prove that the above coloring of edges is a coarsening of $\Theta^*_{S(G)}$. Let $a$ be a vertex of $G$ and $b,c$ its neighbors. Since $G$ is 2-connected, there exists a $b,c$-path $P$ that does not cross $a$. Pick $P$ such that it is shortest possible. Then since $G$ is chordal, $a$ is adjacent to every vertex on $P$, otherwise there exists a shorter path. Denote $P=p_0p_1\ldots p_k$, where $p_0=b$ and $p_k=c$. Then $\{\overline{p_ia}, \overline{p_i}\} \Theta_{S(G)} \{\overline{p_{i+1}a}, \overline{p_{i+1}}\}$, proving that  $\{\overline{ba}, \overline{b}\} \Theta_{S(G)}^* \{\overline{ca}, \overline{c}\}$.

Furthermore, if $ab$ is not an exposed edge in $G$, then $ab$ lies in two maximal cliques. In particular, it lies in two isometrically touching triangles. By Lemma \ref{lem:cycles}, $\{\overline{ab}, \overline{b}\} \Theta_{S(G)}^* \{\overline{a}, \overline{ab}\}$. By transitivity, and the above two facts, all the edges $\{\overline{ab}, \overline{b}\}$, with $a$ being in the same connected component of $G^{-ee}$, are in relation  $\Theta_{S(G)}^*$.

Finally, we prove that no other edge besides the asserted is in $\Theta_{S(G)}^*$. Assume otherwise, and let $\{\bar{a},\overline{ab}\} \Theta_{S(G)}^* \{\bar{c},\overline{cd}\}$ be such that $b$ and $d$ do not lie in the same connected component of $G^{-ee}$. By Lemma \ref{lem:local}, we can assume that $\{\bar{a},\overline{ab}\} \Phi_{S(G)} \{\bar{c},\overline{cd}\}$. But then the edges lie on a 6-cycle, implying that $b=d$. This cannot be.
\qed

\section*{Acknowledgments}
We thank Reza Naserasr and Petru Valicov for fixing a bug in a previous version of Lemma~\ref{lem:separating}.

The authors acknowledge the financial support from the Slovenian Research Agency (research core funding No.\ P1-0297, projects J1-9109 and N1-0095). Kolja Knauer was partially supported by ANR grant DISTANCIA: ANR-17-CE40-0015.


\end{document}